\documentclass[11pt]{amsart}
\usepackage{amsfonts}

\usepackage{amsmath,amsthm,amsfonts,amssymb,latexsym,epsfig}

\usepackage{hyperref}

\newtheorem{theorem}{Theorem}[section]

\newtheorem{lemma}{Lemma}[section]

\newtheorem{con}{Conjecture}[section]

\numberwithin{equation}{section}

\theoremstyle{remark}

\begin{document}
\title{On Edwards--Child's Inequality}
\author[Y.-D. Wu]{Yu-Dong Wu$^*$}
\address[Y.-D. Wu]{Department of Mathematics, Zhejiang Xinchang High School, Shaoxing, Zhejiang 312500, P. R. China.}
\email{yudong.wu@hotmail.com}
\author[Z.-H. Zhang]{Zhi-Hua Zhang}
\address[Z.-H. Zhang]{Department of Mathematics, Shili Senior High School in Zixing, Chenzhou,
 Hunan 423400, P. R. China.} \email{zxzh1234@163.com}
\author[Z.-G. Wang]{Zhi-Gang Wang}
\address[Z.-G. Wang]{School of Mathematics and Computing Science,
Changsha University of Science and Technology, Changsha, Hunan
410076, P. R. China.} \email{zhigwang@163.com }
\thanks{$^*$Corresponding author, e-mail address: yudong.wu@hotmail.com}
\date{January 1, 2008}

\begin{abstract}
In this paper, by making use of one of Chen's theorems and the
method of mathematical analysis, we refine Edwards--Child's
inequality and solve a conjecture posed by Liu \cite{lbq02}.

\end{abstract}

\maketitle



\vskip.20in \section{Introduction and Main Results} \vskip.10in

For $\triangle ABC$, let $a,b,c$ denote the side-lengths, $A,B,C$
the angles, $p$ the semi-perimeter, $R$ the circumradius and $r$ the
inradius, respectively.

Bottema {\it et al.} \cite[pp.50, Theorem 5.8]{b01} once recorded
the following so-called Edwards-Child's inequality in their
monograph:
\begin{equation}\label{01}
p\geqq 3\sqrt{3}\ r.
\end{equation}
\par In 1999, Liu  \cite{lbq01} considered the refinement of
inequality \eqref{01} and posed the problem: Determine the best
constant $t$ for the inequality
\begin{equation}\label{02}
a+b+c\geqq 6\sqrt{3}\
r+t\cdot\frac{r}{R}\cdot(R-2r)\Longleftrightarrow p\geqq 3\sqrt{3}\
r+\frac{t}{2}\cdot\left(1-\frac{2r}{R}\right)r.
\end{equation}
In the same year, Huang \cite{hbc01} solved the above problem and
obtained the best constant $$t=6.829212418\cdots.$$
\par In 2000, Liu \cite{lbq02} posed a conjecture which sharpened
inequality \eqref{01} as follows:
\begin{con}{\rm(LBQ 67)}
Prove or disprove the inequality:
\begin{equation}\label{03}
\sqrt{3}\ p\geqq 10r-r\left(\frac{2r}{R}\right)^5\Longleftrightarrow
p\geqq 3\sqrt{3}\
r+\frac{\sqrt{3}}{3}\left[1-\left(\frac{2r}{R}\right)^5\right]r.
\end{equation}
\end{con}
The main object of this paper is to prove and refine inequality
\eqref{03}, and we obtain the following results.
\begin{theorem}\label{t01}
If $\lambda\leqq f(x_1)\approx 5.977930729$, then
\begin{equation}\label{04}
\sqrt{3}\ p\geqq 10r-r\left(\frac{2r}{R}\right)^\lambda,
\end{equation}
the equality in \eqref{04} is valid if and only if $a:b:c=1:1:1$
with $\lambda\leqq f(x_1)$ or
$a:b:c=2\left(x_1^2-3\right):\left(x_1^2+3\right):\left(x_1^2+3\right)$
with $\lambda=f(x_1)$, where
\begin{equation}\label{09}
f(x):=\frac{\ln{\left[-\frac{2(x^3-5x^2+15)}{x^2-3}\right]}}{\ln{\left[\frac{24(x^2-3)}{(x^2+3)^2}\right]}}\qquad(3<
x<x_0),
\end{equation}
$x_1 \approx 3.067873979\in(3, x_0)$ is the only positive real root
of the equation $f'(x)=0$, and $x_0\approx 4.113537611$ is the
maximal real root of the following equation:
$$x^3-5x^2+15=0.$$
\end{theorem}
\begin{theorem}\label{t02}
If $k\leqq k_0$, then for any triangle, we have
\begin{equation}\label{05}
p\geqq 3\sqrt{3}\ r+k\left[1-\left(\frac{2r}{R}\right)^5\right]r,
\end{equation}
the equality in \eqref{05} is valid if and only if $a:b:c=1:1:1$
with $k\leqq k_0$ or $a:b:c=2t_1:1:1$ with $k= k_0$, where
$k_0\approx 0.6898369707\in\left(\frac{1}{2},\frac{9}{13}\right)$ is
the root of the equation $p_4(k)=0$, $t_1\approx
0.5194285605\in\left(\frac{1}{2},\frac{3}{5}\right) $ is the root of
the equation $p_2(t)=0$, and
\begin{align}
\begin{split}\label{061}
p_4(k)=5&82076609134674072265625k^{40}\\&-347825698554515838623046875k^{38}\\&
+100074581801891326904296875000k^{36}\\&-18656729921698570251464843750000k^{34}\\&
+2406402165103435516357421875000000k^{32}\\&-223881811253170562744140625000000000k^{30}\\&
+13463175828870323610839843750000000000k^{28}\\&-546554928845186341347265625000000000000k^{26}\\&
+44602844570033253018161875000000000000000k^{24}\\&-5130242398470886015317438950000000000000000k^{22}\\&
+805214383330095009369880969748800000000000000k^{20}\\&+8091804003905867976154901735852032000000000000k^{18}\\&
-262076670696960781271382283372767764480000000000k^{16}\\&+483089080872113925346827615868453187944448000000000k^{14}\\&
+76133654804831682593043073614469922630706135040000000k^{12}\\&+4656778296665388277933286028092539858012767959121920000k^{10}\\&
+144475786697302680271016689740018636544036558347226316800k^8\\&+2305862812061518538327375046497265994061106944720616030208k^6\\&
+14131109130840787698067340120866948829204788093469111353344k^4\\&-27802797644590317388762803455393575731200000000000000000000k^2\\&
+9774552621457470122500000000000000000000000000000000000000,
\end{split}
\end{align}
and
\begin{align}\begin{split}\label{062}
p_{2}(t)=&91750400t^{20}-857210880t^{19}+3560046592t^{18}-8616673280t^{17}\\&+13364330496t^{16}
-13704527872t^{15}+9174425600t^{14}-3719282688t^{13}\\&+689182720t^{12}+53528576t^{11}
-36429312t^{10}+533760t^9+207104t^8\\&+75520t^7+32352t^6+11376t^5-320t^4
-104t^3-30t^2-7t-1.
\end{split}\end{align}
\end{theorem}


\vskip.10in
\section{Preliminary Results} \vskip.10in

In order to prove our main results, we shall require the following
three lemmas.
\begin{lemma}{\rm(see \cite{sl01,sl02,yd01})}\label{le01}
The homogeneous inequality $p\geqq (>)f_{1}(R,r)$ in triangle holds
if and only if
\begin{equation}\label{ydeq21}
2(1+t)\sqrt{1-t^2}\geqq (>)f_{1}(1,2t(1-t)),
\end{equation}
where $t=\cos{B}=\cos{C}\in \left[\frac{1}{2},1\right)$. And the
homogeneous inequality $p\leqq (<)f_{1}(R,r)$ holds if and only if
\eqref{ydeq21} is reversed, where $t=\cos{B}=\cos{C}\in
\left(0,\frac{1}{2}\right]$.
\end{lemma}

\begin{lemma}{\rm(see \cite{yhz01,yzh01})}\label{le02}
Given a polynomial $f(x)$ with real coefficients
\begin{equation*}
f(x)=a_{0}x^{n}+a_{1}x^{n-1}+\cdots+a_{n}\nonumber,
\end{equation*}
if the number of the sign changes of the revised sign list of its
discriminant sequence
\begin{equation*}
\left\{D_{1}(f),D_{2}(f),\cdots,D_{n}(f)\right\}\nonumber
\end{equation*}
is $v$, then, the number of the pairs of distinct conjugate
imaginary roots of $f(x)$ equals $v$. Furthermore, if the number of
non-vanishing members of the revised sign list is $l$, then, the
number of the distinct real roots of $f(x)$ equals $l-2v$.
\end{lemma}

\begin{lemma}{\rm(see \cite{yzh01})}\label{le04}
Denote $$ F(x)=a_{0}x^{n}+a_{1}x^{n-1}+\cdots+a_{n},$$ and $$
G(x)=b_{0}x^{m}+b_{1}x^{m-1}+\cdots+b_{m}.\nonumber $$ If $a_{0}\neq
0$ or $b_{0}\neq 0$, then the polynomials $F(x)$ and $G(x)$ have the
common roots if and only if
\begin{equation*}
R(F,G)=
\begin{vmatrix}
  a_{0} & a_{1} & a_{2} & \cdots & a_{n} & 0 & \cdots & 0 \\
  0 & a_{0} & a_{1} & \cdots & a_{n-1} & a_{n} & \cdots & \cdots \\
  \cdots & \cdots & \cdots & \cdots & \cdots & \cdots & \cdots & \cdots \\
  0 & 0 & \cdots & a_{0} & \cdots & \cdots & \cdots & a_{n} \\
  b_{0} & b_{1} & b_{2} & \cdots & \cdots & \cdots & \cdots & 0 \\
  0 & b_{0} & b_{1} & \cdots & \cdots & \cdots & \cdots & 0 \\
  \cdots & \cdots & \cdots & \cdots & \cdots & \cdots & \cdots & \cdots \\
  0 & 0 & 0 & \cdots & b_{0} & b_{1} & \cdots & b_{m} \\
\end{vmatrix}
=0,
\end{equation*}
where $R(F,G)$ is Sylvester's Resultant of $F(x)$ and $G(x)$.
\end{lemma}


\vskip.20in  \section{The Proof of Theorem \ref{t01}} \vskip.10in

\begin{proof}
By Lemma \ref{le01}, we know that inequality \eqref{04} can be
written as follows:
\begin{equation}\label{06}
2\sqrt{3}(1+t)\sqrt{1-t^2}\geqq
20t(1-t)-2t(1-t)[4t(1-t)]^\lambda\qquad \left(\frac{1}{2}\leqq
t<1\right).
\end{equation}

If we let $$x=\sqrt{\frac{3+3t}{1-t}},$$ then
$$t=\frac{x^2-3}{x^2+3}\qquad (x\geqq 3),$$ and inequality \eqref{06} is
equivalent to
\begin{equation}\label{07}
\left[\frac{24(x^2-3)}{(x^2+3)^2}\right]^\lambda\geqq
-\frac{2(x^3-5x^2+15)}{x^2-3}\qquad (x\geqq 3).
\end{equation}

We now split it into three cases to prove.

(i) When $x=3$, inequality \eqref{07} is an identity for $\lambda\in
\mathbb{R}$ and $\triangle{ABC}$ is equilateral.

(ii) Clearly, if $-\frac{2(x^3-5x^2+15)}{x^2-3}\leqq 0$, or $x\geqq
x_{0}\approx 4.113537611$, then inequality \eqref{07} holds for
$\lambda\in \mathbb{R}$.

(iii) When $3< x< x_0$, we have
$$0<\frac{24(x^2-3)}{(x^2+3)^2}< 1,$$
and inequality \eqref{07} is equivalent to
\begin{equation*}\label{08}
\lambda \leqq f(x),
\end{equation*}
where $f(x)$ is given by \eqref{09}. Calculating the derivative for
$f(x)$, we get
\begin{equation*}\label{10}
f'(x)=\frac{x(x-3)(x+3)g(x)}{(x^2+3)(x^2-3)(x^3-5x^2+15)\ln^2{\left[\frac{24(x^2-3)}{(x^2+3)^2}\right]}}.
\end{equation*}where
\begin{equation*}\label{12}
g(x):=x(x^2+3)\ln{\left[\frac{24(x^2-3)}{(x^2+3)^2}\right]}+2(x^3-5x^2+15)\ln{\left[-\frac{2(x^3-5x^2+15)}{x^2-3}\right]}\end{equation*}$$(3<
x<x_0).$$

Let $f'(x)=0$, we obtain
\begin{equation}\label{11}
g(x)=0.
\end{equation}

Now, we prove equation \eqref{11} has only one real root on interval
$(3, x_0)$.

For $g(x)$, we know that
\begin{equation*}\label{13}
g'(x)=3(x^2+1)\ln{\left[\frac{24(x^2-3)}{(x^2+3)^2}\right]}+2x(3x-10)\ln{\left[-\frac{2(x^3-5x^2+15)}{x^2-3}\right]},
\end{equation*}
\begin{align*}
g''(x)=&\frac{2x(x-3)(x+3)(5x^4+6x^3-60x^2-45)}{(x^3-5x^2+15)(x^2+3)(x^2-3)}+6x\ln{\left[\frac{24(x^2-3)}{(x^2+3)^2}\right]}\\&
+(12x-20)\ln{\left[-\frac{2(x^3-5x^2+15)}{x^2-3}\right]},
\end{align*}
\begin{align*}
g'''(x)=&\frac{2}{(x^3-5x^2+15)^2(x^2+3)^2(x^2-3)^2}\big(20x^{13}-113x^{12}-180x^{11}+1500x^{10}\\&
-540x^9-2349x^8+1620x^7-14256x^6-29160x^5+164025x^4+58320x^3\\&-364500x^2-54675\big)+6\ln{\left[\frac{24(x^2-3)}{(x^2+3)^2}\right]}
+12\ln{\left[-\frac{2(x^3-5x^2+15)}{x^2-3}\right]},
\end{align*}
and
\begin{align*}
g^{(4)}(x)=&\frac{4xp(x)}{(x^3-5x^2+15)^3(x^2+3)^3(x^2-3)^3},
\end{align*}
where
\begin{align*}
p(x)=&5x^{18}-78x^{17}+375x^{16}-306x^{15}-705x^{14}-4239x^{13}
+13230x^{12}-16443x^{11}\\&+213435x^{10}-1047330x^9+1356750x^8+4312764x^7-12141495x^6-2460375x^5\\&
+22766670x^4+4100625x^3-32805000x^2-14762250x+61509375.
\end{align*}

The revised sign list of the discriminant sequence of $p(x)$ is
\begin{equation}\label{20}
[1, 1, -1, 1, 1, -1, -1, -1, 1, 1, -1, -1, -1, 1, 1, -1, 1, -1].
\end{equation}
So the number of the sign changes of the revised sign list of
\eqref{20} is $9$. By Lemma \ref{le02}, we know that the number of
the pairs of distinct conjugate imaginary roots of $p(x)$ are $9$.
Namely, $p(x)=0$ has no real root, then $p(x)>0$ for $x\in R$.

 With
$p(x)>0$ for $x\in R$, $g^{(4)}(x)<0$ for $x\in(3,x_0)$, then
$g'''(x)$ is strictly decreasing on $(3,x_0)$. And for
$$g'''(3)=3>0\quad {\rm and}\quad g'''(4)=6\ln{\left(\frac
{312}{361}\right)}-\frac {106269606}{61009}+
12\ln{\left(\frac{2}{13}\right)}<0,$$ we conclude that $g'''(x)=0$
has only one real root $x_4\approx 3.016763142$ on $(3,x_0)$. So
$g''(x)$ is strictly increasing on $(3,x_4)$, and decreasing on
$(x_4,x_0)$.

Similarly, we deduce that $$g^{(j)}(x)=0\qquad (0\leqq j\leqq 2)$$
have only one real root $x_{j+1}\ (0\leqq j\leqq 2)$ on $(3,x_0)$,
respectively.

Therefore, we know that $f'(x)=0$ has only one real root $x_1\approx
3.067873979\in(3,x_0)$, and
$$\lambda_{\max}=f(x)_{\min}=f(x_1)\approx 5.977930729.$$

By Lemma \ref{le01}, we find that the equality in \eqref{04} is
valid when $\triangle{ABC}$ is an isoceles triangle, and if we let
$b=c=1$, then
$$a=2\cos{B}=2\cos{C}=2t_1=\frac{2(x_1^2-3)}{x_1^2+3}.$$ Namely, the
equality in \eqref{04} is valid if and only if
$$a:b:c=2\left(x_1^2-3\right):\left(x_1^2+3\right):\left(x_1^2+3\right)$$ with $\lambda= f(x_1)$.
Furthermore, considering the proof above, we can easily obtain the
required results of the equality in \eqref{04}.

\end{proof}

\vskip.20in  \section{The Proof of Theorem \ref{t02}} \vskip.10in

\begin{proof}
By Lemma \ref{le01}, we know that inequality \eqref{05} is
equivalent to
\begin{equation}\label{41}
2(1+t)\sqrt{1-t^2}\geqq 6\sqrt{3}t(1-t)+k\cdot
2t(1-t)\left\{1-[4t(1-t)]^5\right\}\qquad\left(\frac{1}{2}\leqq
t<1\right).
\end{equation}
For $t=\frac{1}{2}$, \eqref{41} obviously holds and $\triangle{ABC}$
is equilateral. If $\frac{1}{2}<t<1$, then \eqref{41} is equivalent
to
\begin{equation}\label{42}
k\leqq
h(t):=\frac{(1+t)\sqrt{1-t^2}-3\sqrt{3}t(1-t)}{t(1-t)\left\{1-[4t(1-t)]^5\right\}}\qquad
\left(\frac{1}{2}<t<1\right).
\end{equation}
Therefore, we get
\begin{align*}
h'(t)=\big[&5120t^{12}-14336t^{11}+1024t^{10}+30720t^9-25600t^8-10240t^7+19456t^6
-6144t^5\\&+t+1+15360\sqrt{3}\sqrt{1-t^2}(1-t)^5t^6\big]/\big[t^2(1-t)(2t-1)^3(256t^8
-1024t^7\\&+1472t^6-832t^5+80t^4+32t^3+12t^2+4t+1)^2\sqrt{1-t^2}\
\big].
\end{align*}

Let $h'(t)=0$, we obtain
\begin{align}\label{43}\begin{split}
5120t^{12}&-14336t^{11}+1024t^{10}+30720t^9-25600t^8-10240t^7+19456t^6\\&
-6144t^5+t+1+15360\sqrt{3}\sqrt{1-t^2}(1-t)^5t^6=0.\end{split}
\end{align}
It is easy to see that the roots of equation \eqref{43} must be the
roots of the following equation:
\begin{equation}\label{44}
p_{2}(t)(t+1)(2t-1)^3=0,
\end{equation}
where $ p_2(t)$ is given by \eqref{062}.

The revised sign list of the discriminant sequence of $p_{2}(t)$ is
\begin{equation}\label{45}
[1, 1, 1, -1, -1, 1, -1, -1, -1, 1, 1, 1, -1, 1, 1, 1, -1, -1, 1,
1].
\end{equation}
So the number of the sign changes of the revised sign list of
\eqref{45} is $8$, by Lemma \ref{le02}, the equation $p_{2}(t)=0$
has $4$ distinct real roots. By using the function "fsolve()" in
Maple 9.0, we can obtain the approximation of $2$ distinct real
roots $t_1\approx 0.5194285605 $ and $t_2\approx 0.8281776966$ on
open interval $\left(\frac{1}{2},1\right)$. But $t_2$ does not
satisfy equation \eqref{43}, which implies that it  is an extraneous
root. Hence, it follows that
\begin{equation}\label{451}
h(t)_{\min}=h(t_1)=\frac{(1+t_1)\sqrt{1-t_1^2}-3\sqrt{3}t_1(1-t_1)}{t_1(1-t_1)\{1-[4t_1(1-t_1)]^5\}}\approx
0.6898369707\in \left(\frac{1}{2},\frac{9}{13}\right).
\end{equation}

We now prove $h(t_1)$ is the root of equation \eqref{06}. We
consider the nonlinear algebraic equation system as follows:
\begin{equation}\label{46}
\begin{cases}
p_{2}(t_1)=0,\\
u^2+t_1^2-1=0,\\
v^2-27t_1^2(1-t_1)^2=0,\\
(1+t_1)u-v-\left\{t_1(1-t_1)\left\{1-[4t_1(1-t_1)]^5\right\}\right\}k=0,\\
\end{cases}
\end{equation}
or
\begin{equation}\label{461}
\begin{cases}
p_{2}(t_1)=0,\ \ \ \ \\
p_3(t_1,k)=0,
\end{cases}
\end{equation}
where $ p_2(t)$ is given by \eqref{062} and
\begin{align*}\label{1111}
p_3(t,k)=&68719476736k^4t^{42}-1374389534720k^4t^{41}\\&
+13022340841472k^4t^{40}-77687368450048k^4t^{39}+327083234426880k^4t^{38}\\&
-1032364109070336k^4t^{37}+2532355667394560k^4t^{36}-4938039864328192k^4t^{35}\\&
+7763189089435648k^4t^{34}-9918303552143360k^4t^{33}+10328492122046464k^4t^{32}\\&
-8752260867686400k^4t^{31}+5995600248045568k^4t^{30}-3278892823478272k^4t^{29}\\&
+1402218741760000k^4t^{28}-453762372075520k^4t^{27}+105664190873600k^4t^{26}\\&
-16867120906240k^4t^{25}+2458252738560k^4t^{24}-944229580800k^4t^{23}\\&
+65536k^2(-52+7344385k^2)t^{22}-1441792k^2(-24+119075k^2)t^{21}\\&
+4587520k^2(8441k^2-34)t^{20}-65536k^2(65905k^2-6238)t^{19}\\&
+20480k^2(-33404+4115k^2)t^{18}-32768k^2(2431k^2-22973)t^{17}\\&
+2048k^2(29895k^2-262244)t^{16}-2048k^2(12617k^2-113240)t^{15}\\&
+256k^2(-196132+23549k^2)t^{14}-2560k^2(39k^2-1232)t^{13}\\&
-512k^2(1929+71k^2)t^{12}-256k^2(-1576+55k^2)t^{11}-16k^2(109k^2-2660)t^{10}\\&
-64k^2(51k^2-1274)t^9+40k^2(-296+7k^2)t^8+8k^2(-788+11k^2)t^7\\&
+k^2(25k^2-3076)t^6+2k^2(-400+3k^2)t^5+k^2(k^2-166)t^4\\&-20k^2t^3+(49-2k^2)t^2+14t+1.
\end{align*}
It's easy to see that $h(t_{1})$ is also the solution of the
nonlinear algebraic equation system \eqref{461}. By Lemma
\ref{le04}, we get
\begin{equation}
R(p_2,p_3)=mp_4(k)p_5(k)=0,
\end{equation}
where
\begin{align*}
m=&57698221920758648437318785957934716092387171383190921139891196583\\&
60844636195941737072193545378747970974059915467496354369852962661\\&
24545032425816659873003224168956415237423809191615340664942700611\\&
579699989665614710761572731854962712529967513600000000,
\end{align*}
$ p_4(k)$ is given by \eqref{061}, and
\begin{align}
\begin{split}
&p_5(k)=8860655573197291232645511627197265625k^{40}\\&
-16912080685346637403103522956371307373046875k^{38}\\&
+5435741841027337227212034165859222412109375000k^{36}\\&
-1022146637691011812806702821850776672363281250000k^{34}\\&
+110744571814857445274894881859313964843750000000000k^{32}\\&
-9369389460811282973529757765591658325195312500000000k^{30}\\&
+282441031609128684811039782683995118457031250000000000k^{28}\\&
-10380425845408883215082367863858860726041875000000000000k^{26}\\&
+73340535725547628521708894311600819307265000000000000000k^{24}\\&
-265177139233592420219133919392417723053323929920000000000000k^{22}\\&
+11791633182359761811369142031608045030413197651840000000000000k^{20}\\&
+190075444754710412043307126977386168137413366964111360000000000k^{18}\\&
-109270095333365249190090664013722276751701441874134340403200000000k^{16}\\&
-7750407940399239639127577820761730582290929415800425022750720000000k^{14}\\&
-128016712123970225804454667626572657168697762262885566533284659200000k^{12}\\&
+2859604104165956733160109025625828353920705420748548328012916981760000k^{10}\\&
+110267783320883416462258578945595318603200831554078609670205562316390400k^8\\&
+745324132711516724121423152749040579524897407250894331064162322304991232k^6\\&
-7375168578676180534548979250039227170865907645609839475147982859847860224k^4\\&
-80207674110325900841453359343883145581345268086466150400000000000000000000k^2\\&
+395739228446715359134735466442056102500000000000000000000000000000000000000.
\end{split}\end{align}
The revised sign list of the discriminant sequence of $p_{5}(k)$ is
\begin{align}\label{47}\begin{split}
[1, 1, 1, 1, -1, 1, 1, -1, 1, &1, 1, 1, -1, 1, -1, 1, 1, 1, 1, -1,
1, 1, 1, -1, 1, 1, 1, 1, \\&1, -1, -1, -1, 1, 1, -1, 1, 1, 1, 1,
1].\end{split}
\end{align}
So the number of the sign changes of the revised sign list of
\eqref{47} is $16$, by Lemma \ref{le02}, the equation $p_{5}(k)=0$
has $8$ distinct real roots. And making use of the function
"realroot()"\cite{dm01} in Maple 9.0, we know that the equation
$p_{5}(k)=0$ has $8$ distinct real roots in the following intervals:
\begin{equation*}\label{48}
\left[2, \frac{5}{2}\right],\ \ \left[\frac{5}{2}, 3\right], \ \ [4,
8], \ \ [1024, 2048], \ \ \left[-\frac{5}{2}, -2\right],\ \
\left[-3, -\frac{5}{2}\right],\ \ [-8, -4],\ \ [-2048, -1024].
\end{equation*}
Therefore, the equation $p_{5}(k)=0$ has no real root on interval
$\left(\frac{1}{2},\frac{9}{13}\right)$.\par From \eqref{451}, we
find that $h(t_1)$ is the root of the equation $p_{4}(k)=0$. Namely,
$h(t_1)$ is the root of the equation \eqref{061}.

Furthermore, by the similar arguments of the last part of the proof
of Theorem \ref{t01}, we can easily get the necessary and sufficient
conditions on the equality in \eqref{05}.

We thus complete the proof of Theorem \ref{t02}.
\end{proof}

\vskip.10in

\begin{center}{\sc Acknowledgements}
\end{center}
\vskip.10in

The present investigation was supported, in part, by the {\it Hunan
Provincial Natural Science \linebreak Foundation} under Grant
05JJ30013 and the {\it Scientific Research Fund of Hunan Provincial
\linebreak Education Department} under Grant 05C266 of People's
Republic of China.

\end{document}